\documentclass[12pt]{article}
\usepackage{amssymb}

\usepackage{amsmath,latexsym,amsfonts}
\usepackage{graphicx}
\usepackage{srcltx}
\usepackage{amscd}

\newtheorem{lemma}{Lemma}[section]
\newtheorem{coro}[lemma]{Corollary}

\newtheorem{thm}[lemma]{Theorem}
\newtheorem{defn}[lemma]{Definition}

\makeatletter\@addtoreset{equation}{section}
\renewcommand\theequation{\thesection.\@arabic\c@equation}

\begin{document}
\begin{center}
{\LARGE Liouville type theorems for transversally harmonic and biharmonic maps}

 \renewcommand{\thefootnote}{}
\renewcommand{\thefootnote}{\arabic{footnote}}
\setcounter{footnote}{0}

\vspace{0.5 cm} {\large Min Joo Jung and Seoung Dal Jung}
\end{center}
\vspace{0.5cm}
\noindent {\bf Abstract.} In this paper, we study the Liouville type theorems for transversally harmonic and biharmonic maps on foliated Riemannian manifolds.


\section{Introduction}


Let $(M,\mathcal F)$ and $(M',\mathcal F')$  be
foliated Riemannian manifolds and let $\phi:M\to M'$ be a smooth  foliated map,
i.e., $\phi$ is a smooth leaf-preserving map.
Then $\phi$ is said to be {\it
transversally harmonic}  if the transversal tension field $\tau_b(\phi)={\rm tr}_Q{\tilde\nabla}_{\rm tr}d_T\phi$
vanishes, where $d_T\phi=d\phi|_Q$ and $Q$ is the normal bundle of $\mathcal F$ (see [\ref{JJ2},\ref{KW1},\ref{KW2}] for details). When $\mathcal F$ is minimal, a transversally harmonic map is  a critical point of the transversal energy $E_B(\phi)$ [\ref{JJ2}], which is given by
\begin{align*}
E_B(\phi) =\frac12\int_M |d_T\phi|^2 \mu_M,
\end{align*}
where $\mu_M$ denotes the volume form on $M$. If $\mathcal F$ is not minimal, a transversally harmonic is not a critical point of $E_B(\phi)$. In fact, S. Dragomir and  A. Tommasoli [\ref{DT}] called such maps as $(\mathcal F,\mathcal F')$-harmonic maps, i.e., a critical point of the transversal energy. Trivially, two definitions are equivalent when $\mathcal F$ is minimal. 
The smooth map $\phi$ is said to be {\it transversally biharmonic} if the transversal bitension field $(\tau_2)_b(\phi)=J_\phi^T(\tau_b(\phi)$ vanishes, where $J_\phi^T$ is the generalized Jacobi operator along $\phi$ (see 
[\ref{CW},\ref{JU1}] for details). If $\mathcal F$ is minimal, then a trasversally biharmonic map is a critical point of the transversal bienergy $E_2(\phi)$, where
\begin{align*}
E_2(\phi)=\frac12\int_M |\tau_b(\phi)|^2\mu_M.
\end{align*}
Transversally harmonic and biharmonic maps are generalizations of  harmonic and biharmonic maps because transversally harmonic and biharmonic maps are just harmonic and biharmonic maps on the point foliation, respectively. For more information about transversally harmonic and biharmonic maps, see [\ref{CW},\ref{JJ2},\ref{JU1},\ref{KW1},\ref{KW2}].  For harmonic maps, the classical Liouville theorem is well-known. Namely, any bounded harmonic function defined on the whole plane must be constant. The classical Liouville theorem has been improved in several cases [\ref{JU},\ref{SY},\ref{YA}]. In this article, we study the Liouville type theorems for the transversally harmonic and biharmonic map.  Now we consider the following conditions on $(M,g,\mathcal F)$ and $(M',g',\mathcal F')$.

\bigskip
(C1) All 
leaves of $\mathcal F$ are compact and  the mean curvature form $\kappa$ of $\mathcal F$ is  
bounded, coclosed.

(C2) The transversal sectional curvature of $\mathcal F'$ is nonpositive. 

\bigskip

\noindent Then we have the following Liouville type theorem on a foliated Riemannian manifold.

\bigskip
\noindent{\bf Theorem A.} {\it
Let $(M,g,\mathcal F)$ be a complete foliated Riemannian manifold with {\rm Vol}$(M)=\infty$ satisfying  $(C1)$ and let $(M',g',\mathcal F')$ be a foliated Riemannian manifold satisfying $(C2)$.  Assume that    the the transversal Ricci curvature of $\mathcal F$ is nonnegative. Then any transversally harmonic map $\phi:M\to M'$ of $E_B(\phi) <\infty$  is transversally constant, i.e., the induced map between leaf spaces is constant.}

\bigskip
\noindent Note that any transversally harmonic map is transversally biharmonic. But the converse does not hold. In fact, S. D. Jung [\ref{JU1}] proved that on a compact foliated manifold, the converse holds under some condition. For transversally biharmonic map on a complete foliated Riemannian manifold, we have the following theorem.

\bigskip
\noindent
{\bf Theorem B.} {\it 
Let $(M,g,\mathcal F)$ be a complete foliated Riemannian manifold with {\rm Vol}$(M)=\infty$ satisfying  $(C1)$ and let $(M',g',\mathcal F')$ be a foliated Riemannian manifold satisfying $(C2)$.  

\noindent$(1)$ Every transversally biharmonic map $\phi:M\to M'$ of $E_2(\phi) <\infty$ is transversally harmonic.

\noindent$(2)$ If the transversal Ricci curvature of $\mathcal F$ is nonnegative,  then every transversally biharmonic map $\phi:M\to M'$ of $E_B(\phi)+E_2(\phi) <\infty$ is transversally constant.   }

\bigskip
\noindent
When $\mathcal F$ is a point foliation, Theorem A and Theorem B have been found in [\ref{SY}] and [\ref{BF}], respectively.
\section{Preliminaries}
Let $(M,g,\mathcal F)$ be a $(p+q)$-dimensional Riemannian
manifold with a foliation $\mathcal F$ of codimension $q$ and a complete
bundle-like metric $g$ with respect to $\mathcal F$. 
Let $TM$ be the tangent bundle of $M$, $T\mathcal F$ its integrable
subbundle given by $\mathcal F$, and  $Q=TM/T\mathcal F$ the corresponding
normal bundle of $\mathcal F$.
 Then we have an exact sequence of vector bundles
\begin{align}\label{eq1-1}
 0 \longrightarrow T\mathcal F \longrightarrow
TM_{\buildrel \longleftarrow \over \sigma }^{\buildrel \pi \over
\longrightarrow} Q \longrightarrow 0,
\end{align}
where $\pi:TM\to Q$ is a projection and $\sigma:Q\to T\mathcal F^\perp$ is a bundle map satisfying
$\pi\circ\sigma=id$.  
    Let $g_Q$ be the holonomy invariant metric
on $Q$ induced by $g$, i.e., $L_Xg_Q=0$ for any  vector field $X\in T\mathcal F$, where
$L_X$ is the transverse Lie derivative [\ref{KT1}]. Let $R^Q$, $K^Q$ and ${\rm Ric}^Q$  be the transversal curvature tensor, transversal sectional curvature and transversal Ricci operator
of $\mathcal F$ with respect to the transversal Levi-Civita connection $\nabla^Q\equiv \nabla$ in $Q$ [\ref{TO}], respectively.
A differential form $\omega\in \Omega^r(M)$ is {\it basic} if $
i(X)\omega=0$ and $i(X)d\omega=0$ for all $X\in T\mathcal F$. 
Let $\Omega_B^r(\mathcal F)$ be the set of all basic $r$-forms on
$M$. Then $\Omega^r(M)=\Omega_B^r(\mathcal F)\oplus \Omega_B^r(\mathcal F)^\perp$ [\ref{LO}].
  Now, we recall the star operator $\bar *:\Omega_B^r (\mathcal F)\to \Omega_B^{q-r}(\mathcal F)$ given by [\ref{KT2},\ref{PR}]
\begin{align}
\bar *\omega= (-1)^{p(q-r)}*(\omega\wedge\chi_{\mathcal F}),\quad\forall \omega\in\Omega_B^r(\mathcal F),
\end{align}
where $\chi_{\mathcal F}$ is the characteristic form of $\mathcal F$ and $*$ is the Hodge star operator associated to $g$.
For any basic forms $\omega,\theta \in \Omega_B^r(\mathcal F)$, it is well-known [\ref{PR}] that $\omega\wedge\bar *\theta = \theta\wedge\bar *\omega$ and $\bar *^2\omega = (-1)^{r(q-r)}\omega$.  
Let $\nu$ be the transversal volume form, i.e., $*\nu =\chi_{\mathcal F}$ and  $\langle\cdot,\cdot\rangle$ be the pointwise inner product on $\Omega_B^r(\mathcal F)$, which is given by
\begin{align}
\langle\omega,\theta\rangle \nu = \omega\wedge\bar *\theta
\end{align}
for any basic forms $\omega,\theta \in \Omega_B^r(\mathcal F)$. 
 Trivially  $\mu_M=\nu\wedge\chi_{\mathcal F}$ is the volume form with respect to $g$.  Now, let  the operator $d_B$ be the restriction of $d$ to the basic forms, i.e., $d_B= d|_{\Omega_B^*(\mathcal F)}$.  It is well-known that on complete foliated Riemannian manifolds,  $d_B\kappa_B=0$ [\ref{PJ}]. 
 Let  $d_t = d_B-\kappa_B\wedge $, $\delta_t = (-1)^{q(r+1)+1}\bar * d_B\bar *$ and 
\begin{align}
\delta_B\omega = (-1)^{q(r+1)+1}\bar * d_t\bar *\omega = \delta_t\omega + i(\kappa_B^\sharp)\omega,
\end{align}
where  $(\cdot)^\sharp$ is the $g_Q$-dual vector field of $(\cdot)$ and  $\kappa_B$ is the basic part of the mean curvature form of $\mathcal F$.  Then  $\int_M\langle d_B\omega,\theta\rangle\mu_M=\int_M\langle \omega,\delta_B\theta\rangle\mu_M$ for any $\omega\in\Omega_{B,o}^r(\mathcal F)$ or $\theta\in\Omega_{B,o}^{r+1}(\mathcal F)$, where $\Omega_{B,o}^*(\mathcal F)$ is the subspace of $\Omega_B^*(\mathcal F)$ composed of forms with compact support.  Generally,  $\delta_B$ is not a restriction of $\delta$ on $\Omega_B^r (\mathcal F)$, i.e., $\delta_B \ne \delta|_{\Omega_B^r (\mathcal F)}$, where $\delta$ is the formal adjoint of $d$. But $\delta_B\omega=\delta\omega$ for any basic $1$-form $\omega$. Hence $\Delta^M|_{\Omega_B^0(\mathcal F) } =\Delta_B$ [\ref{KT2}], where $\Delta^M$ is the Laplacian on $M$ and $\Delta_B$ is the basic Laplacian acting on $\Omega_B^r (\mathcal F)$ which is given by
\begin{align}
\Delta_B =d_B\delta_B +\delta_B d_B.
\end{align}
Let $V(\mathcal F)$ be the space of all transversal infinitesimal automorphisms $Y$ of $\mathcal F$, i.e., $[Y,Z]\in T\mathcal F$ for all $Z\in T\mathcal F$. 
For any  $Y\in V(\mathcal F)$, we define the bundle map $A_Y:\Gamma(\Lambda^r
Q^*)\to\Gamma(\Lambda^r Q^*)$  [\ref{KT1}]
by
\begin{align}\label{eq1-13}
A_Y\omega =L_Y\omega-\nabla_Y\omega.
\end{align}
Then $A_Y$ preserves the basic forms and depends only on $\pi(Y)$. Moreover, for any vector field
$Y\in V(\mathcal F)$, if we define $A_Y :\Gamma Q\to \Gamma Q$ by $\omega(A_Ys)=-(A_Y\omega)(s)$ for any $\omega\in \Gamma Q^*$, then
$A_Ys = L_Ys-\nabla_Ys$. Since $L_Ys=\pi[Y,Y_s]$ for $Y_s =\sigma(s)\in T\mathcal F^\perp$ [\ref{KT1}], 
\begin{equation}
A_Y s = -\nabla_{Y_s}\pi(Y).
\end{equation} 
Let
 $\{E_a\}(a=1,\cdots,q)$ be a local orthonormal basic  frame of $Q$ and $\theta^a$ a $g_Q$-dual $1$-form to $E_a$. We define  $\nabla_{\rm tr}^*\nabla_{\rm tr}:\Omega_B^r(\mathcal F)\to \Omega_B^r(\mathcal F)$ by
\begin{align}\label{eq1-12}
\nabla_{\rm tr}^*\nabla_{\rm tr} =-\sum_a \nabla^2_{E_a,E_a}
+\nabla_{\kappa_B^\sharp},
\end{align}
where $\nabla^2_{X,Y}=\nabla_X\nabla_Y -\nabla_{\nabla^M_XY}$ for
any $X,Y\in TM$ and $\nabla^M$ is the Levi-Civita connection with respect to $g$.  The operator $\nabla_{\rm tr}^*\nabla_{\rm tr}$
is positive definite and formally self adjoint on $\Omega_{B,o}^r(\mathcal F)$ [\ref{Jung}]. 
Then the  generalized Weitzenb\"ock type formula on $\Omega_B^r(\mathcal F)$ is given by [\ref{Jung}]
\begin{align}\label{2-7}
\Delta_B \omega = \nabla_{\rm tr}^*\nabla_{\rm tr} \omega
+ F(\omega) + A_{\kappa_{B}^\sharp} \omega
\end{align}
for any $\omega\in\Omega_B^r(\mathcal F)$,
where  $F=\sum_{a,b=1}^{q}\theta^a\wedge i(E_b)  R^Q(E_b,E_a)$.
For any basic-harmonic form $\omega\in \Omega_B^*(\mathcal F)$, i.e., $\Delta_B\omega=0$, we have [\ref{Jung}] that
\begin{equation}
-\frac12\Delta_B |\omega|^2 =|\nabla_{\rm tr}\omega|^2 + \langle A_{\kappa_B^\sharp}\omega,\omega\rangle +\langle
F(\omega),\omega\rangle.
\end{equation}

\section{Generalized maximum principle}
Let $(M,g,\mathcal F)$ be a complete foliated Riemannian manifold, i.e., manifold with a Riemannian  foliation $\mathcal F$ and a complete bundle-like metric $g$ with respect to $\mathcal F$. 
Now, we consider a smooth function $\mu$ on $\mathbb R$ satisfying
\begin{align*}
(i) \ 0\leq \mu(t)\leq 1\ {\rm on}\ \mathbb R,\quad
(ii)\ \mu(t)=1\ \ {\rm for}\ t\leq 1,\quad
(iii)\ \mu(t)=0\ \ {\rm for}\ t\geq 2.
\end{align*}
Let $x_0$ be a point in $M$. For each point $y \in M$, we
denote by $\rho(y)$ the distance between leaves through $x_0$ and $y$. For any real number $l>0$, we define a Lipschitz continuous function $\omega_l$  on $M$ by  
\begin{align*}
\omega_l(y)=\mu(\rho(y)/l).
\end{align*} 
Trivially, $\omega_l$ is a basic function. 
   Let $B(l) = \{ y \in M | \rho(y) \leq l \}$. Then $\omega_l$ satisfies
the following properties:
\begin{center}
 $\left.
 \begin{array}{ll}
0 \leq \omega_l(y) \leq 1 & \text{for any} \ y \in M\\
\text{supp}\ \omega_l \subset B(2l) & \\
\omega_l(y)=1 & \text{for any} \  y \in B(l)\\
\lim_{l \rightarrow \infty} \omega_l = 1 & \\
|d\omega_l | \leq \frac{C}{l} & \text{almost everywhere on}\ M, 
\end{array} \right.$
\end{center}
where $C$ is a positive constant independent of $l$ [\ref{YO1}]. Hence $\omega_l\psi$ has compact support for any basic form $\psi\in\Omega_B^*(\mathcal F)$ and $\omega_l\psi \to \psi$ (strongly) when $l\to \infty$.

Note that for any basic function $f$, $\Delta^Mf =\Delta_B f$ [\ref{KT2}]. Hence we have the following theorems.
\begin{thm} Let $(M,g,\mathcal F)$ be a complete foliated Riemannian manifold. If a basic function  $f$ is basic-subharmonic, i.e., $\Delta_B f\leq 0$, with $\int_M |df|<\infty$, then $f$ is basic-harmonic.
\end{thm}
{\bf Proof.} This follows from the result in [\ref{YA2}, p 660]. $\Box$

\begin{thm} Let $(M,g,\mathcal F)$ be a complete foliated Riemannian manifold. If a nonnegative basic function $f$ is basic-subharmonic, i.e., $\Delta_B f\leq 0$, with $\int_M f^p <\infty\ (p>1)$, then $f$ is constant.
\end{thm}
{\bf Proof.} This follows from Theorem 3 in [\ref{YA2}, p 663]. $\Box$

   Now we prove the generalized maximum principle.
\begin{thm} 
Let $(M,g,\mathcal F)$ be a complete foliated Riemannian manifold whose all leaves are compact. Assume that $\kappa_B$ is bounded and coclosed.   Then a nonnegative basic function $f$ such that $(\Delta_B-\kappa_B^\sharp)f\leq 0$ with $\int_M f^p <\infty\ (p>1)$ is constant.
\end{thm}
{\bf Proof.} Let $u=f^{p/2}$.  By a direct calculation, we have
\begin{align}\label{3-5-1}
u(\Delta_B-\kappa_B^\sharp)u ={p\over 2}f^{p-1}(\Delta_B-\kappa_B^\sharp)f-{p(p-2)\over 4}f^{p-2}|d_B f|^2.
\end{align}
By the assumption, we have
\begin{align}\label{3-5-2}
u(\Delta_B-\kappa_B^\sharp)u\leq -{p-2\over p}|d_B u|^2.
\end{align}
On the other hand,  we have
\begin{align}\label{3-5-3}
\int_{B(2l)}\langle\omega_l^2 u,\Delta_Bu\rangle &=2\int_{B(2l)}\langle\omega_l d_B u, ud_B\omega_l\rangle +\int_{B(2l)} | \omega_l d_B u|^2.
\end{align}
So, from (\ref{3-5-2}) and (\ref{3-5-3}), we have that
\begin{align}\label{3-5-4}
{2(p-1)\over p}\int_{B(2l)}|\omega_l d_Bu|^2 \leq -2\int_{B(2l)} \langle\omega_l d_B u,u d_B\omega_l\rangle +\int_{B(2l)}\langle \omega_l^2 u,\kappa_B^\sharp(u)\rangle.
\end{align}
From (\ref{3-5-4}) and the Schwarz's inequality, we have that for any real number $\epsilon>0$,
\begin{align*}
{p-1\over p}\int_{B(2l)}|\omega_l d_Bu|^2&\leq \int_{B(2l)}|\langle\omega_l d_Bu,u d_B\omega_l\rangle|+ \frac12\int_{B(2l)}\langle \omega_l^2 u,\kappa_B^\sharp(u)\rangle \\
&\leq {\epsilon\over 2}\int_{B(2l)}|\omega_l d_B u|^2 +{1\over 2\epsilon}\int_{B(2l)}| ud_B\omega_l|^2+ \frac12\int_{B(2l)}\langle \omega_l^2 u,\kappa_B^\sharp(u)\rangle \\
&\leq  {\epsilon\over 2}\int_{B(2l)}| \omega_l d_B u|^2+{C^2\over 2\epsilon l^2}\int_{B(2l)} u^2+ \frac12\int_{B(2l)}\langle \omega_l^2 u,\kappa_B^\sharp(u)\rangle.
\end{align*}
Hence we have
\begin{align}
({p-1\over p}-{\epsilon\over 2})\int_{B(2l)}|\omega_l d_Bu|^2
&\leq  {C^2\over 2\epsilon l^2}\int_{B(2l)} u^2+ \frac12\int_{B(2l)}\langle \omega_l^2 u,\kappa_B^\sharp(u)\rangle.
\end{align}
On the other hand, since $\langle\omega_l^2 u,\kappa_B^\sharp(u)\rangle = {1\over 2}\{\kappa_B^\sharp(\omega_l^2 u^2) -2\langle d\omega_l,(\omega_l u^2 )\kappa_B\rangle\}$, we have  that from the assumptions of $\kappa_B$, i.e., $|\kappa_B|<\infty$ and $\delta_B\kappa_B=0$,  
\begin{align*}
\int_{B(2l)}|\langle \omega_l^2 u,\kappa_B^\sharp(u)\rangle|&\leq\int_{B(2l)}|d\omega_l| |(\omega_l u^2 )\kappa_B |\\
&\leq {C\over l}{\max} (|\kappa_B|)\int_{B(2l)} \omega_l u^2.
\end{align*}
Since $\int_M u^2 =\int_M f^p<\infty$, if we let $l\to \infty$, then
\begin{align}
 \int_{B(2l)}\langle \omega_l^2 u,\kappa_B^\sharp(u)\rangle \to 0.
\end{align}
From (3.5) and (3.6),  if we let $l\to \infty$, then 
\begin{align}\label{3-6-1}
({p-1\over p}-{\epsilon\over 2})\int_M| d_B u|^2 \leq 0.
\end{align}
If we choose $0<\epsilon< {2(p-1)\over p}$, then  from (\ref{3-6-1}),
\begin{align*}
\int_M | d_B u|^2 =0.
\end{align*}
Hence $d_B u=0$ and so $d_Bf=0$, i.e., $f$ is constant. $\Box$ 

\bigskip
\noindent{\bf Remark.}  On a compact foliated Riemannian manifold, Theorem 3.3 was proved in [\ref{JLR}].

\section{The proof of Theorem A} 
Let $(M,  g,\mathcal F)$  and $(M', g',\mathcal F')$  be
two Riemannian manifolds with foliations $\mathcal F$ and $\mathcal F'$, respectively.
  Let $\nabla$ and $\nabla'$ be the transverse
Levi-Civita connections on $Q$ and $Q'$, respectively.  Let $\phi:(M,g,\mathcal
F)\to (M', g',\mathcal F')$ be a smooth  foliated map,
i.e., $\phi$ is a smooth leaf-preserving map. Equivalently, $d\phi(T\mathcal F)\subset T\mathcal F'$.  We define $d_T\phi:Q \to Q'$  by
\begin{align}\label{4-1}
d_T\phi := \pi' \circ d \phi \circ \sigma.
\end{align}
Then $d_T\phi$ is a section in $ Q^*\otimes
\phi^{-1}Q'$, where $\phi^{-1}Q'$ is the pull-back bundle on $M$.  Let $\nabla^\phi$
and $\tilde \nabla$ be the connections on $\phi^{-1}Q'$ and
$Q^*\otimes \phi^{-1}Q'$, respectively. 
The {\it
transversal tension field} of $\phi$ is defined by
\begin{align}\label{4-3}
\tau_b(\phi)={\rm tr}_{Q}\tilde \nabla d_T
\phi=\sum_{a=1}^{q}(\tilde\nabla_{E_a} d_T\phi)(E_a),
\end{align}
where $\{E_a\}(a=1,\cdots,q)$ is a local orthonormal basic frame on $Q$.
Trivially, the transversal tension field $\tau_b(\phi)$ is a
section of $\phi^{-1}Q'$.
A foliated map $\phi: (M, g,\mathcal F) \to (M', g',\mathcal F')$ is said to be 
{\it transversally harmonic}  if the transversal tension field
vanishes, i.e., $\tau_b(\phi)=0$ [\ref{JJ2}].  And the {\it transversal energy} of $\phi$ on a compact domain $\Omega$ is defined by
\begin{align}
E_B(\phi;\Omega)=\frac12\int_\Omega |d_T\phi|^2\mu_M,
\end{align}
where $|d_T\phi|^2 =\sum_a g_{Q'}(d_T\phi(E_a),d_T\phi(E_a))\in \Omega_B^0(\mathcal F)\ [\ref{KG}]$. Then we have the  first variational formula [\ref{JJ2}]
\begin{align}
{d\over dt} E_B(\phi_t;\Omega)|_{t=0}=-\int_\Omega \langle V,\tau_b(\phi)-d_T\phi(\kappa_B^\sharp)\rangle\mu_M,
\end{align}
where $V={d\phi_t \over dt} |_{t=0}$ is the normal variation vector field with a foliated variation $\{\phi_t\}$ of $\phi$. Hence if $\mathcal F$ is minimal, then the transversal harmonic map is a critical point of the transversal energy $E_B(\phi;\Omega)$ of $\phi$ supported in a compact domain $\Omega$. 
Let $\Omega_B^*(E)\equiv \Omega_B^*(\mathcal F)\otimes E$, where  $E\equiv\phi^{-1}Q'$.
Then we define $d_\nabla :\Omega_B^r(E)\to \Omega_B^{r+1}(E)$ by
\begin{align}\label{2-6}
d_\nabla( \omega\otimes s)=  d_B \omega\otimes s+ (-1)^r\omega\wedge \nabla^\phi_{\rm tr} s
\end{align}
for any $\omega\in\Omega_B^r(\mathcal F)$ and $s\in \Gamma E$. 
Let $\delta_\nabla$ be the formal adjoint of $d_\nabla$ on $\Omega_{B,o}^*(E)$, the space of the compact supports. We define the Laplacian $\Delta$ on $\Omega_B^r(E)$  by
\begin{align}\label{2-6-1}
  \Delta=d_\nabla \delta_\nabla +\delta_\nabla d_\nabla.
  \end{align}
 The operator $A_Y$ is extended to $\Omega_B^r(E)$ [\ref{JJ2}]. That is, for any $\omega\otimes s \in \Omega_B^r (E)$,
 $A_Y(\omega\otimes s) = A_Y\omega \otimes s$. 
Then  we have the generalized Weitzenb\"ock type formula.  

\begin{thm} $[\ref{JJ2}]$ Let $\phi:(M,g,\mathcal F)\to (M',g',\mathcal F')$ be a smooth foliated map. Then
\begin{align*}
\frac12\Delta_B|d_T\phi|^2=\langle\Delta d_T\phi,d_T\phi\rangle-|\tilde\nabla_{\rm tr}d_T\phi|^2-\langle A_{\kappa_B^\sharp}d_T\phi,d_T\phi\rangle-\langle F(d_T\phi),d_T\phi\rangle,
\end{align*}
where 
\begin{align}\label{3-7}
\langle F(d_T\phi),d_T\phi\rangle &=\sum_a g_{Q'}(d_T\phi({\rm Ric}^Q(E_a)),d_T\phi(E_a))\notag\\
&-\sum_{a,b}g_{Q'}(R^{Q'}(d_T\phi(E_a),d_T\phi(E_b))d_T\phi(E_b),d_T\phi(E_a)).
\end{align}
\end{thm}
Note that  for a smooth foliated map $\phi:(M,g,\mathcal F)\to (M',g',\mathcal F')$ [\ref{JJ2}],
\begin{align}\label{3-5}
d_\nabla d_T\phi=0,\quad \delta_\nabla d_T\phi=-\tau_b(\phi)+i(\kappa_B^\sharp)d_T\phi.
\end{align}
Hence we have the following corollary.
\begin{coro} $[\ref{JJ2}]$ Let $\phi:(M,g,\mathcal F)\to (M',g',\mathcal F')$ be a transversally harmonic map. 
Then 
\begin{align}\label{eq4-13}
\frac12 (\Delta_B-\kappa_B^\sharp) |d_T\phi|^2 =-|\tilde\nabla_{\rm tr}d_T\phi|^2 -\langle F(d_T\phi),d_T\phi\rangle.
\end{align}
\end{coro}
{\bf The proof of Theorem A.} Note that $\frac12\Delta_B |d_T\phi|^2 = |d_T\phi| \Delta_B |d_T\phi| - |d_B |d_T \phi||^2$. From  Corollary 4.2, we have
\begin{align*}
|d_T\phi|(\Delta_B-\kappa_B^\sharp) |d_T\phi| &= |d_B |d_T\phi||^2 -|\tilde \nabla_{\rm tr} d_T\phi |^2-\langle F(d_T\phi),d_T\phi\rangle.
\end{align*}
Since $|\tilde\nabla_{\rm tr} d_T\phi | \geq |d_B |d_T\phi||$ (Kato's inequality [\ref{BE}]),  we have
\begin{align*}
|d_T\phi|(\Delta_B-\kappa_B^\sharp) |d_T\phi| \leq -\langle F(d_T\phi),d_T\phi\rangle.
\end{align*}
By the assumptions of the curvatures, $\langle F(d_T\phi),d_T\phi\rangle\geq 0$, which means 
 $ (\Delta_B-\kappa_B^\sharp) |d_T\phi|\leq 0$. Hence, by  Theorem 3.3,  $|d_T\phi|$ is constant. Since ${\rm Vol}(M)$ is  infinite and $E_B(\phi)<\infty$,  we have $d_T\phi=0$. Hence $\phi$ is transversally constant. $\Box$ 


\section {The proof of Theorem B}
 Let $\phi:(M,g,\mathcal
F)\to (M', g',\mathcal F')$ be a smooth  foliated map. 
 The {\it transversal bitension field $(\tau_2)_b(\phi)$} of $\phi$ is defined by
\begin{align}\label{6-1}
(\tau_2)_b(\phi)=J_\phi^T(\tau_b(\phi)),
\end{align}
 where the generalized Jacobi operator $J_\phi^T: \phi^{-1}Q'\to  \phi^{-1}Q'$ {\rm along $\phi$ is defined by
\begin{align}\label{4-6}
J^T_\phi(s)=(\nabla_{\rm tr}^\phi)^*(\nabla_{\rm tr}^\phi)s-\nabla_{\kappa_B^\sharp}^\phi s-{\rm tr}_Q R^{Q'}(s,d_T\phi)d_T\phi
\end{align}
for any $s\in \phi^{-1}(Q')$ [\ref{JU1}].
\begin{defn} {\rm [\ref{CW}]} {\rm Let $\phi:(M,g,\mathcal F)\to (M',g',\mathcal F')$ be a smooth foliated map. Then $\phi$ is said to be} transversally biharmonic {\rm if the transversal bitension field vanishes, i.e., $(\tau_2)_b(\phi)=0$.}
\end{defn}
Trivially, $\phi$ is a transversally biharmonic map if and only if the transversal tension field $\tau_b(\phi)$ is a generalized Jacobi field along $\phi$. 

\bigskip
\noindent
{\bf The proof of Theorem B.}
Note that for any $s\in\phi^{-1}Q'$, we have
 \begin{align}\label{4-8}
 \frac12\Delta_B |s|^2 = \langle (\nabla_{\rm tr}^\phi)^* (\nabla_{\rm tr}^\phi) s,s\rangle-|\nabla_{\rm tr}^\phi s|^2.
 \end{align} 
Since $\phi:M\to M'$ is transversally biharmonic, from (5.1) and (5.2),
 \begin{align}\label{5-3}
(\nabla_{\rm tr}^\phi)^*(\nabla_{\rm tr}^\phi)\tau_b(\phi)-\nabla_{\kappa_B^\sharp}^\phi \tau_b(\phi)-{\rm tr}_Q R^{Q'}(\tau_b(\phi),d_T\phi)d_T\phi=0.
\end{align}
From (\ref{4-8}) and (\ref{5-3}), we have
\begin{align}
\frac12(\Delta_B -\kappa_B^\sharp)|\tau_b(\phi)|^2 = \langle{\rm tr}_Q R^{Q'}(\tau_b(\phi),d_T\phi)d_T\phi,\tau_b(\phi)\rangle - |\nabla_{\rm tr}^\phi \tau_b(\phi)|^2.
\end{align} 
Since $\frac12\Delta_B f^2 = f\Delta_B f - |d_B f|^2$ for any basic function $f$, Eq. (5.5) implies that
\begin{align*}
|\tau_b(\phi)| (\Delta_B -\kappa_B^\sharp)|\tau_b(\phi)| =&\langle{\rm tr}_Q R^{Q'}(\tau_b(\phi),d_T\phi)d_T\phi,\tau_b(\phi)\rangle\\
 &+ |d_B |\tau_b(\phi)||^2- |\nabla_{\rm tr}^\phi \tau_b(\phi)|^2.
\end{align*}
By assumption of the transversal sectional curvature of $\mathcal F'$, i.e., $K^{Q'}\leq 0$ and the Kato's inequality [\ref{BE}], we have
\begin{align}
(\Delta_B -\kappa_B^\sharp)|\tau_b(\phi)|\leq 0.
\end{align}
Since $\int_M |\tau_b(\phi)|^2<\infty$,  by Theorem 3.3,  $|\tau_b(\phi)|$ is constant. Hence ${\rm Vol}(M)=\infty$ implies that $\tau_b(\phi)=0$, i.e., $\phi$ is transversally harmonic. This complete the proof of (1).
For the proof of (2), it is trivial from Theorem A and Theorem B (1). $\Box$

\bigskip
\noindent{\bf Acknowledgements.}  \noindent{\bf Acknowledgements.}  This research was supported by  the National Research Foundation of Korea(NRF) grant funded
       by the Korea government (MSIP) (NRF-2015R1A2A2A01003491).

\noindent Department of Mathematics and Research Institute for Basic Sciences, Jeju National University, Jeju 690-756, Korea

\noindent {\it E-mail address} : niver486@jejunu.ac.kr

\noindent {\it E-mail address} : sdjung@jejunu.ac.kr


\begin{thebibliography}{[9]}

\bibitem{Lop}\label{LO} J. A. Alvarez {L\'opez},
\emph{The basic component of the mean curvature of Riemannian
foliations}, Ann. Global Anal. Geom. 10 (1992), 179-194.

\bibitem{BF} \label{BF} P. Baird, A. Fardoun and S. Ouakkas, \emph{Liouville-type theorems for biharmonic maps between Rieannian manifolds}, Adv. Calc. Var. 3 (2010), 49-68.



\bibitem{be}\label{BE} P. B\'erard, \emph{A note on Bochner type theorems for complete manifolds}, Manuscripta Math. 69 (1990), 261-266.


\bibitem{cw}\label{CW} Y.-J. Chiang and R. Wolak, \emph{Transversally biharmonic maps between foliated Riemannian manifolds}, International J. Math. 19 (2008), 981-996.

\bibitem{DT}\label{DT} S. Dragomir and A. Tommasoli, \emph{Harmonic maps of foliated Riemannian manifolds}, Geom. Dedicata 162 (2013), 191-229.

\bibitem{jjj}\label{JJ2} M. J. Jung and S. D. Jung, \emph{On transversally harmonic maps of foliated Riemannian manifolds}, J. Korean Math. Soc. 49 (2012), 977-991.


\bibitem{j}\label{JU} S. D. Jung, \emph{Harmonic maps of complete Riemannian manifolds}, Nihonkai Math. J. 8 (1997), 147-154.

\bibitem{Jung}\label{Jung} S. D. Jung,
\emph{The first eigenvalue of the transversal Dirac operator}, J.
Geom. Phys. 39 (2001), 253--264.

\bibitem{Jung}\label{JU1} S. D. Jung, \emph{Variation formulas for transversally harmonic and biharmonic maps}, J. Geom. Phys. 70 (2013), 9-20.

\bibitem{JLR} \label{JLR} S. D. Jung, K. R. Lee and K. Richardson, \emph{Generalized Obata theorem and its applications on foliations}, J. Math. Anal. Appl. 376 (2011), 129-135.

\bibitem{k}\label{KG} A. E. Kacimi Alaoui and E. Gallego Gomez, \emph{Applications harmoniques feuillet\'ees}, Illinois J. Math. 40 (1996), 115-122.

\bibitem{Kamber2}\label{KT1} F. W. Kamber and Ph.
Tondeur, \emph{Infinitesimal automorphisms and second variation of
the energy for harmonic foliations}, T\^ohoku Math. J. 34 (1982),
525-538.

\bibitem{jj}\label{KT2} F. W. Kamber and Ph. Tondeur, {\it De Rham-Hodge theory for Riemannian foliations}, Math. Ann. 277 (1987), 415-431.

\bibitem{konderak}\label{KW1} J. Konderak and R. Wolak, \emph{On transversally harmonic maps between manifolds with Riemannian foliations}, Quart. J. Math. Oxford Ser.(2) 54 (2003), 335-354.

\bibitem{konderak1}\label{KW2} J. Konderak and R. Wolak, \emph{Some remarks on transversally harmonic maps}, Glasgow Math. J. 50 (2008), 1-16.


\bibitem{PJ}\label{PJ} J. S. Pak and S. D. Jung, \emph{A transversal Dirac operator and some vainshing theorems on a complete foliated Riemannian manifold}, Math. J. Toyama Univ. 16 (1993), 97-108.

\bibitem{PR}\label{PR} E. Park and K. Richardson, \emph{The basic Laplacian of a Riemannian foliation}, Amer. J. Math. 118 (1996), 1249-1275.




\bibitem{S}\label{SY} R. Schoen and S. T. Yau, {\it Harmonic maps and the topology of stable hypersurfaces and manifolds of nonnegative Ricci curvature,} Comm. Math. Helv. 51 (1976), 333-341.



\bibitem{Tond1}\label{TO} Ph. Tondeur,
\emph{Geometry of foliations}, Basel: Birkh\"auser Verlag, 1997.


\bibitem{Yau}\label{YA} S. T. Yau, {\it Harmonic functions on complete Riemannian manifolds}, Comm. Pure Appl. Math. 28 (1975), 201-228.

\bibitem{Yau2}\label{YA2} S. T. Yau, {\it Some function-theoretic properties of complete Riemannian manifold and their applications to geometry}, Indiana Univ. Math. J. 25 (1976), 659-670.


\bibitem{Yorozu}\label{YO1} S. Yorozu, \emph{Notes on suare-integrable cohomology spaces on certain foliated manifolds}, Trans. Amer. Math. Soc. 255 (1979), 329-341.




\end{thebibliography}
\end{document}